# Alternative proofs of ten inequalities

A.N. Lepes

Department of Mathematics,
The Republican Specialized Physics-Mathematics Secondary Boarding School
named after O. Zhautykov,
36 Bukhar Zhyrau Str., 050000, Almaty, Kazakhstan
Tel.: (727) 249-95-87, Fax: (727) 395-01-77
e-mail: adilsultan01@gmail.com

This article offers different proofs of ten inequalities from those already published. So that the readers can see for themselves, the tasks specified in the condition of the source and classical inequalities which used in previously published proofs.

Source of inspiration for the creation of new proof of inequalities for the author was the article of K. Li (see. [2]).

*Example* 1 (Junior Balkan Mathematical Olympiad 2012, problem 1). Let *a*, *b*, and *c* be positive real numbers such that $a+b+c=1$. Prove that

$$\frac{a}{b}+\frac{a}{c}+\frac{b}{c}+\frac{b}{a}+\frac{c}{b}+\frac{c}{a}+6 \geq 2\sqrt{2}\left(\sqrt{\frac{1-a}{a}}+\sqrt{\frac{1-b}{b}}+\sqrt{\frac{1-c}{c}}\right).$$

*Solution.* Let $f(x) = \frac{1-x}{x} - 2\sqrt{\frac{2(1-x)}{x}}$, where $x \in (0,1)$. Note that if $a = b = c = \frac{1}{3}$, our inequality becomes an equality. Next construct the tangent equation at point $x_0 = \frac{1}{3}$:

$$y = f\left(\frac{1}{3}\right) + f'\left(\frac{1}{3}\right)\left(x - \frac{1}{3}\right) = -2 - 0 \cdot \left(x - \frac{1}{3}\right) = -2.$$

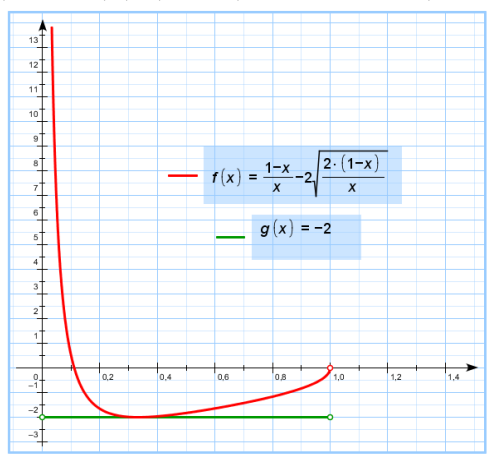

Fig. 1

We now prove that (see Figure 1)

$$\frac{1-x}{x} - 2\sqrt{\frac{2(1-x)}{x}} \geq -2 \tag{1}$$

is true for any $x \in (0, 1)$. Indeed, we can write inequality (1) as follows

$$\frac{1+x}{x} \geq \frac{2\sqrt{2(1-x)}}{\sqrt{x}}. \qquad (2)$$

After multiplying both sides of (2) by $x^2$, squaring, expanding and combining like terms, we have

$$(3x-1)^2 \geq 0,$$

which is clearly true for any $x\in(0, 1)$. Next, using (1), we obtain

$$\frac{1-b}{b}+\frac{1-c}{c}+\frac{1-a}{a} - 2\sqrt{2}\left(\sqrt{\frac{1-a}{a}}+\sqrt{\frac{1-b}{b}}+\sqrt{\frac{1-c}{c}}\right) \geq -6.$$

<div align="right">QED■</div>

The use of a tangent to prove the inequality is complicated by the need to choose a point of tangency. Naturally, the selection of points of tangency of the tangent type, and in accordance with the inequality, requires proof. For example, if one wants on the right-hand side an expression of the form $A+k(a+b+c+d)$, one needs a point of tangency $x_0$ so that $4(f(x_0) - f'(x_0)x_0) = A$. We prove this by the following example.

*Example* 2 (China West Mathematics Invitation, 2004, problem 3, Xiong Bin and Lee Peng Yee (2007): **187**, problem 3). Find all real numbers $k$, such that the inequality
$$a^3 + b^3 + c^3 + d^3 + 1 \geq k(a+b+c+d)$$
holds for any $a, b, c, d\in[-1; +\infty)$.

*Solution.* Let $S=a+b+c+d$, $f(x)=x^3$, where $x\geq-1$. We choose $x_0$ such that $4(f(x_0) - f'(x_0)x_0) = -1$. That is, $x_0=0.5$. Let us construct the tangent equation at this point:

$$y = f(0.5) + f'(0.5)(x-0.5) = \frac{1}{8}+\frac{3}{4}\left(x-\frac{1}{2}\right) = \frac{3x-1}{4}.$$

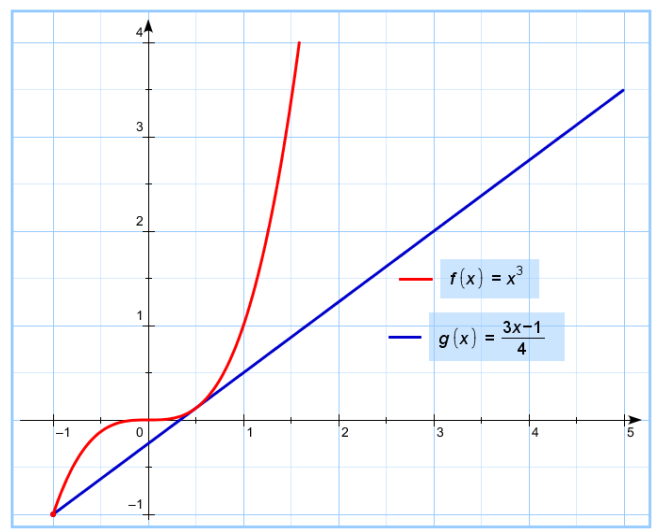

Fig. 2

Note that the inequality (see Figure 2)

$$x^3 \geq \frac{3x-1}{4} \qquad (3)$$

is equivalent to $(x-0.5)^2(x+1) \geq 0$. Therefore, (3) is true for any $x \geq -1$. Applying this inequality, we have

$$a^3 + b^3 + c^3 + d^3 \geq -1 + \frac{3}{4}(a+b+c+d).$$

Thus, $k=0.75$ solves the problem. However, $k$ cannot be greater than 0.75, for when $a=b=c=d=0.5$, the inequality

$$a^3 + b^3 + c^3 + d^3 + 1 \geq k(a+b+c+d)$$

would not hold, contradicting the choice of $k$.

*Example* 3 (Chetkovski, 2012: **185**, problem 25). Let $a, b, c$ be positive real numbers. Prove the inequality

$$\sqrt{\frac{a}{b+c}} + \sqrt{\frac{b}{c+a}} + \sqrt{\frac{c}{a+b}} > 2.$$

*Solution*. Assume that $a+b+c=3$. In this case, our inequality can be written as

$$\sqrt{\frac{a}{3-a}} + \sqrt{\frac{b}{3-b}} + \sqrt{\frac{c}{3-c}} > 2.$$

Let $f(x) = \sqrt{\frac{x}{3-x}}$, with $x \in (0,3)$. Note that if $a=b=1.5$, $c=0$, and our inequality becomes an equality. Hence, we need to set up the equation of the tangent so that it passes through two points on the graph of $f$ with abscissas 0 and 1.5; the graph of the function $f$ would not fall below the tangent, of course, if such tangent exists. Construct the tangent equation to the graph $f$ at the point $x_0 = 1.5$:

$$y = f(1.5) + f'(1.5)(x-1) = 1 + \frac{2}{3}(x-1.5) = \frac{2}{3}x.$$

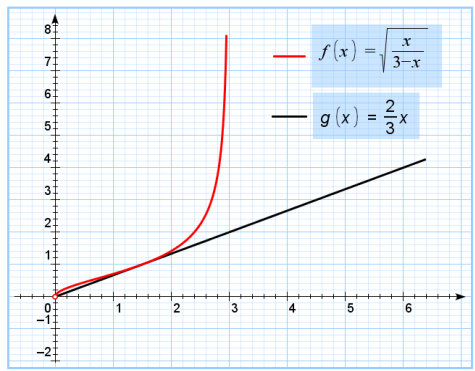

Fig. 3

The inequality $f(x) \geq g(x)$ is true, because (see Figure 3)

$$\sqrt{\frac{x}{3-x}} \geq \frac{2}{3}x \Leftrightarrow x \cdot \left(\frac{4}{9}x^2 - \frac{4}{3}x + 1\right) \geq 0 \Leftrightarrow x \cdot \left(\frac{2}{3}x - 1\right)^2 \geq 0,$$

which holds for any $x \in (0, 3)$. Then applying it thrice with numbers $a$, $b$, and $c$, we find that

$$\sqrt{\frac{a}{3-a}} + \sqrt{\frac{b}{3-b}} + \sqrt{\frac{c}{3-c}} \geq \frac{2}{3}a + \frac{2}{3}b + \frac{2}{3}c = \frac{2}{3}(a+b+c) = 2.$$

Because $a, b, c$ are positive, the equality is possible only if $a = b = c = \frac{3}{2}$. However, $\frac{3}{2} + \frac{3}{2} + \frac{3}{2} \neq 3$ and the equality is impossible. Consequently,

$$\sqrt{\frac{a}{3-a}} + \sqrt{\frac{b}{3-b}} + \sqrt{\frac{c}{3-c}} > 2.$$

QED∎

A tangent might be not a straight line but instead a degree function, for example.

*Example* 4 (Chetkovski, 2012: **189**, problem 63). Let $a$, $b$, $c$ and $d$ be four positive real numbers such that $a^2 + b^2 + c^2 + d^2 = 1$. Prove the inequality

$$\sqrt{1-a} + \sqrt{1-b} + \sqrt{1-c} + \sqrt{1-d} \geq \sqrt{a} + \sqrt{b} + \sqrt{c} + \sqrt{d}.$$

*Solution.* Let $f(x) = \sqrt{1-x} - \sqrt{x}$, $g(x) = kx^2 + m$, where $0 < x < 1$. Numbers $k$ and $m$ are chosen such that $f\left(\frac{1}{2}\right) = g\left(\frac{1}{2}\right)$ and $f'\left(\frac{1}{2}\right) = g'\left(\frac{1}{2}\right)$. That is, $k$ and $m$ satisfy the following equations $0 = \frac{k}{4} + m$ and $-\sqrt{2} = k$. Hence, $g(x) = -\sqrt{2}\left(x^2 - \frac{1}{4}\right)$. We prove that for every $x \in (0,1)$ we have (see Figure 4)

$$\sqrt{1-x} - \sqrt{x} \geq -\sqrt{2}\left(x^2 - \frac{1}{4}\right) \quad (4)$$

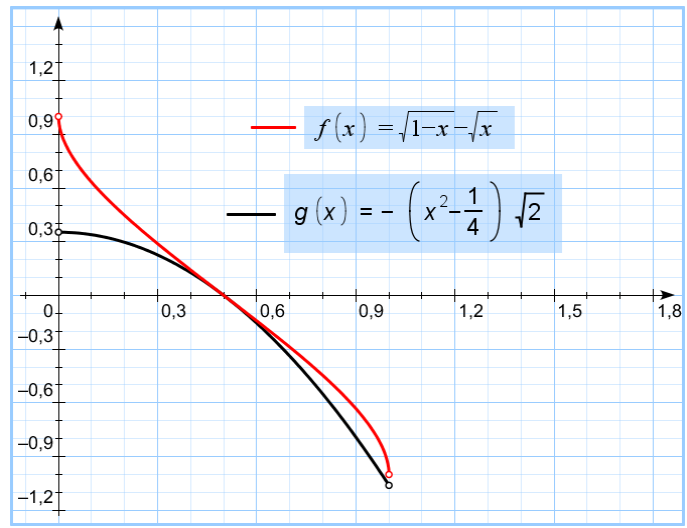

Fig. 4

After performing some equivalent transformations of inequality (4):

$$\frac{1-2x}{\sqrt{1-x}+\sqrt{x}} \geq \frac{(1-2x)(1+2x)}{2\sqrt{2}} \Leftrightarrow \frac{(1-2x)\left(2\sqrt{2}-(1+2x)(\sqrt{1-x}+\sqrt{x})\right)}{2\sqrt{2}(\sqrt{1-x}+\sqrt{x})} \geq 0 \Leftrightarrow$$

$$\Leftrightarrow \frac{(1-2x)\left(2\sqrt{2}-(1+2x)\left(\frac{1-2x}{\sqrt{1-x}+\sqrt{x}}+2\sqrt{x}\right)\right)}{2\sqrt{2}(\sqrt{1-x}+\sqrt{x})} \geq 0 \Leftrightarrow$$

$$\Leftrightarrow \frac{(1-2x)\left(\frac{2x-1}{\sqrt{1-x}+\sqrt{x}}(1+2x)+2(\sqrt{2}-\sqrt{x}-2x\sqrt{x})\right)}{2\sqrt{2}(\sqrt{1-x}+\sqrt{x})} \geq 0 \Leftrightarrow$$

$$\Leftrightarrow \frac{(1-2x)\left(\frac{2x-1}{\sqrt{1-x}+\sqrt{x}}(1+2x)+2(1-\sqrt{2x})(\sqrt{2}x+\sqrt{x}+\sqrt{2})\right)}{2\sqrt{2}(\sqrt{1-x}+\sqrt{x})} \geq 0 \Leftrightarrow$$

$$\Leftrightarrow \frac{(1-\sqrt{2x})^2(1+\sqrt{2x})\left(2(\sqrt{2}x+\sqrt{x}+\sqrt{2})-\frac{(1+\sqrt{2x})(1+2x)}{\sqrt{1-x}+\sqrt{x}}\right)}{2\sqrt{2}(\sqrt{1-x}+\sqrt{x})} \geq 0,$$

$$\Leftrightarrow \frac{(1-\sqrt{2x})^2(1+\sqrt{2x})(2(\sqrt{2}x+\sqrt{x}+\sqrt{2})\sqrt{1-x}+\sqrt{2x}-1)}{2\sqrt{2}(\sqrt{1-x}+\sqrt{x})^2} \geq 0,$$

we find that the last inequality is true for every $x \in (0,1)$ because

$$\sqrt{2}(\sqrt{1-x}+\sqrt{x}) \geq \sqrt{2}\sqrt{1-x+x} = \sqrt{2} > 1.$$

Hence, inequality (4) is valid for the specified $x$. Hence,

$$\sqrt{1-a}+\sqrt{1-b}+\sqrt{1-c}+\sqrt{1-d}-\sqrt{a}-\sqrt{b}-\sqrt{c}-\sqrt{d} \geq$$

$$\geq -\sqrt{2}(a^2+b^2+c^2+d^2-1) = 0.$$

QED■

*Example* 5 (Chetkovski, 2012: **187**, problem 37). Let $a$, $b$, and $c$ be positive real numbers such that $abc \geq 1$. Prove the inequality $\left(a+\frac{1}{a+1}\right)\left(b+\frac{1}{b+1}\right)\left(c+\frac{1}{c+1}\right) \geq \frac{27}{8}$.

*Solution*. Let $f(x) = x + \frac{1}{x+1}$ and $g(x) = kx^m$, where $x>0$. Numbers $k$ and $m$ are chosen such that $f(1) = g(1)$ and $f'(1) = g'(1)$. That is, $k$ and $m$ satisfy the following equations $\frac{3}{2} = k$, $\frac{3}{4} = km$. Hence, $g(x) = \frac{3}{2}\sqrt{x}$. As inequality (see Figure 5)

$$x + \frac{1}{x+1} \geq \frac{3}{2}\sqrt{x} \qquad (5)$$

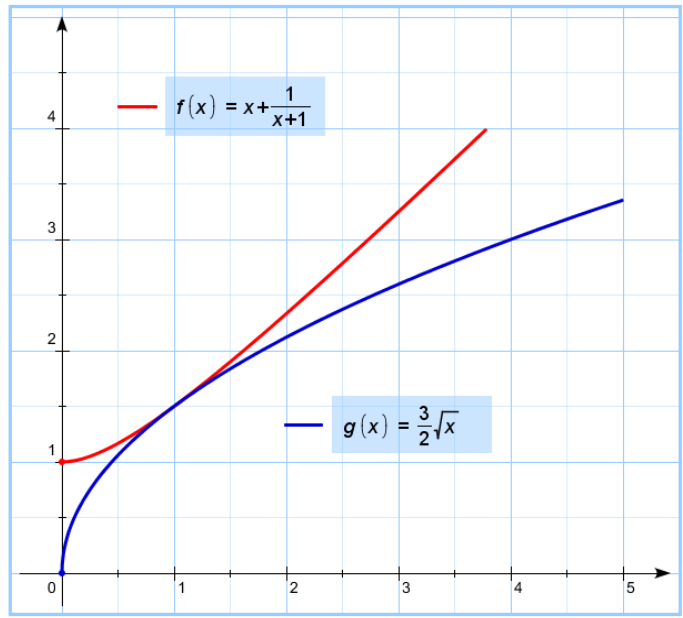

Fig. 5

is equivalent to

$$(\sqrt{x} - 1)^2(2x + \sqrt{x} + 2) \geq 0,$$

which is obviously true for any $x>0$, then (5) is valid for any positive number $x$. Applying (5), then in accordance with the statement of the problem, we have

$$\left(a + \frac{1}{a+1}\right)\left(b + \frac{1}{b+1}\right)\left(c + \frac{1}{c+1}\right) \geq \frac{27}{8}\sqrt{abc} \geq \frac{27}{8}.$$

QED■

As is well-know, one of the applications of inequalities is in finding the largest or smallest values of functions. However, when considering the sum of several different functions, the points of tangency to these functions are chosen such that their derivatives coincide at their respective points. This is to ensure that the relevant coefficients of the variables are the same.

*Example* 6 (Chetkovski, 2012: **39**, Ex.4.5; Cauchy-Schwarz inequality). Let $a$, $b$, and $c$ be positive real numbers. Determine the minimal value of $\frac{3a}{b+c} + \frac{4b}{c+a} + \frac{5c}{a+b}$.

*Solution*. Assume that $a+b+c=1$. Consider the three functions $f_1(x) = 3g(x)$, $f_2(x) = 4g(x)$, and $f_2(x) = 5g(x)$, where $g(x) = \frac{x}{1-x}$, $x\in(0,1)$. Because each $x\in(0,3)$ satisfies the inequality $g''(x) = \frac{3\cdot 2}{(1-x)^3} > 0$, function $g(x)$, and hence the functions $f_1$, $f_2$, $f_3$, are everywhere convex downward on $(0,+\infty)$. Hence, the graphs of these functions are not below any tangent drawn at the point along the positive abscissa. That is, for any $x_1$, $x_2$, $x_3$, $x>0$, we have the following inequalities ($i=1,2, 3$)

$$f_i(x) \geq f_i(x_1) + f_i'(x_i)(x - x_i). \tag{6}$$

Choose positive numbers $x_1$, $x_2$, $x_3$ satisfying the equalities

$$x_1 + x_2 + x_3 = 3,$$

$$f_1'(x_1) = f_2'(x_2) = f_3'(x_3).$$

That is, numbers $x_1$, $x_2$, $x_3$ satisfy the following system of equations (see Figure 6)

$$\begin{cases} x_1 + x_2 + x_3 = 1, \\ \dfrac{\sqrt{3}}{1-x_1} = \dfrac{2}{1-x_2} = \dfrac{\sqrt{5}}{1-x_3}. \end{cases} \tag{7}$$

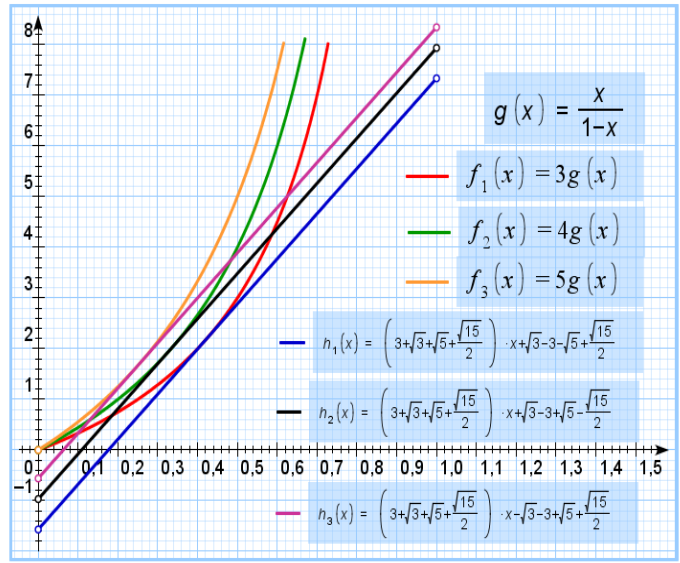

Fig. 6

Let $t = \dfrac{1-x_1}{\sqrt{3}}$, then $x_1 = 1 - \sqrt{3}t$, $x_2 = 1 - 2t$, $x_3 = 1 - \sqrt{5}t$. Hence, from the first equation of (7), we have $t = \dfrac{2}{\sqrt{3}+2+\sqrt{5}}$. Hence, $x_1 = \dfrac{2+\sqrt{5}-\sqrt{3}}{\sqrt{3}+2+\sqrt{5}}$, $x_2 = \dfrac{\sqrt{3}-2+\sqrt{5}}{\sqrt{3}+2+\sqrt{5}}$, $x_3 = \dfrac{\sqrt{3}+2-\sqrt{5}}{\sqrt{3}+2+\sqrt{5}}$. Adding the inequalities of the form (6), for $x$ equal to $x_1$, $x_2$, $x_3$, we have

$$\frac{3a}{1-a} + \frac{4b}{1-b} + \frac{5c}{1-c} = f_1(a) + f_2(b) + f_3(c) \geq$$

$$\geq f_1(x_1) + f_2(x_2) + f_3(x_3) + f_1'(x_1)(a + b + c - x_1 - x_2 - x_3) =$$

$$= \frac{3x_1}{1-x_1} + \frac{4x_2}{1-x_2} + \frac{5x_3}{1-x_3} = \frac{1}{t}\left(\sqrt{3}x_1 + 2x_2 + \sqrt{5}x_3\right) =$$

$$= \frac{\sqrt{3}(2+\sqrt{5}-\sqrt{3}) + 2(\sqrt{3}-2+\sqrt{5}) + \sqrt{5}(\sqrt{3}+2-\sqrt{5})}{2} =$$

$$= \sqrt{15} + 2\sqrt{3} + 2\sqrt{5} - 6 = \frac{\left(\sqrt{3} + 2 + \sqrt{5}\right)^2}{2} - 12.$$

Thus,

$$\min_{a,b,c>0} \left( \frac{3a}{b+c} + \frac{4b}{c+a} + \frac{5c}{a+b} \right) = \frac{\left(\sqrt{3} + 2 + \sqrt{5}\right)^2}{2} - 12 = f(x_1, x_2, x_3).$$

Depending on the complexity of the problem, the proof might require a few applications of the methods or a few well-known inequalities.

*Example* 7 (Chetkovski, 2012: **196**, problem 130). Let $a$, $b$, and $c$ be positive real numbers such that $ab+bc+ca=3$. Prove the inequality

$$(a^7 - a^4 + 3)(b^5 - b^2 + 3)(c^4 - c + 3) \geq 27.$$

Solution. Let $f_1(x) = x^7 - x^4 + 3$, $f_2(x) = x^5 - x^2 + 3$, $f_3(x) = x^4 - x + 3$, $g_1(x) = k_1 x^3 + m_1$, $g_2(x) = k_2 x^3 + m_2$, and $g_3(x) = k_3 x^3 + m_3$, where $x>0$. Note that at $a=b=c=1$, our inequality becomes an equality. Numbers $k_1$, $m_1$, $k_2$, $m_2$, $k_3$, $m_3$ are chosen such that $f_1(1) = g_1(1)$, $f_1'(1) = g_1'(1)$, $f_2(1) = g_2(1)$, $f_2'(1) = g_2'(1)$, $f_3(1) = g_3(1)$, $f_3'(1) = g_3'(1)$. However, $f_1(1) = f_2(1) = f_3(1) = 3$ and $f_1'(1) = f_2'(1) = f_3'(1) = 3$, consequently,

$$g_1(x) = g_2(x) = g_3(x) = x^3 + 2.$$

As inequalities $(x-1)^2(x^2+x+1) \geq 0$ and $(x-1)^2(x^2+x+1)(x+1) \geq 0$ hold for every $x>0$, then the following inequality holds (see Figure 7)

$$x^7 - x^4 + 3 \geq x^5 - x^2 + 3 \geq x^4 - x + 3 \geq x^3 + 2.$$

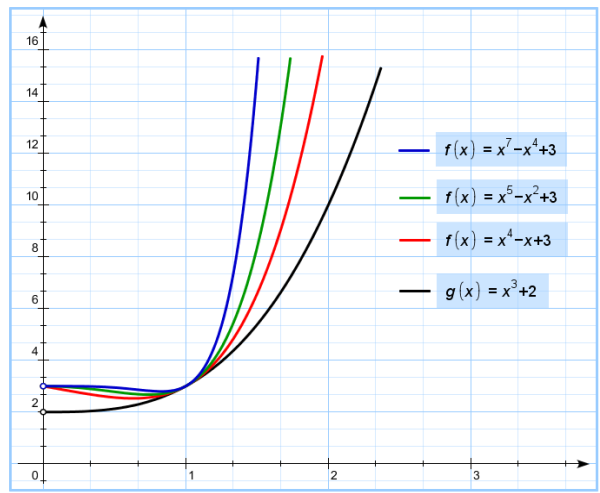

(Fig. 7)

Hence, applying the Hölder inequality and the known inequality

$$(a+b+c)^2 \geq 3(ab+bc+ca),$$

we get

$$(a^7 - a^4 + 3)(b^5 - b^2 + 3)(c^4 - c + 3) \geq (a^3 + 2)(b^3 + 2)(c^3 + 2) =$$

$$= (a^3 + 1^3 + 1^3)(1^3 + b^3 + 1^3)(1^3 + 1^3 + c^3) \geq (a + b + c)^3 \geq$$

$$\geq \left(3(ab + bc + ca)\right)^{\frac{3}{2}} = 27.$$

It should be noted that inequality $(a + b + c)^2 \geq 3(ab + bc + ca)$ is true because it can be represented as

$$(a - b)^2 + (b - c)^2 + (c - a)^2 \geq 0.$$

QED■

*Example* 8 (Chetkovski, 2012: **56**, Ex.5.11). Let $a, b, c, d > 0$ be real numbers. Find the minimum value of the expression

$$\frac{a}{b+c+d} + \frac{b}{c+d+a} + \frac{c}{d+a+b} + \frac{d}{a+b+c} + \frac{b+c+d}{a} + \frac{c+d+a}{b} + \frac{d+a+b}{c} + \frac{a+b+c}{d}.$$

*Solution.* Let

$$f(a,b,c,d) = \frac{a}{b+c+d} + \frac{b}{c+d+a} + \frac{c}{d+a+b} + \frac{d}{a+b+c} + \frac{b+c+d}{a} + \frac{c+d+a}{b} + \frac{d+a+b}{c} + \frac{a+b+c}{d}.$$

Assume that $a+b+c+d=4$. In this case, the function $f$ can be written as

$$f(a,b,c,d) = \frac{a}{4-a} + \frac{b}{4-b} + \frac{c}{4-c} + \frac{d}{4-d} + \frac{4-a}{b} + \frac{4-b}{b} + \frac{4-c}{c} + \frac{4-d}{d}.$$

Consider $g(x) = \frac{x}{4-x} + \frac{4-x}{x}$, $x \in (0,4)$. Construct the tangent equation to the graph $g$ at the point $x_0 = 1$:

$$y = g(1) + g'(1)(x - 1) = \frac{10}{3} - \frac{32}{9}(x - 1) = \frac{62 - 32x}{9}.$$

We prove that for every $x \in (0,4)$ we have (see Figure 8)

$$\frac{x}{4-x} + \frac{4-x}{x} \geq \frac{62-32x}{9}. \tag{8}$$

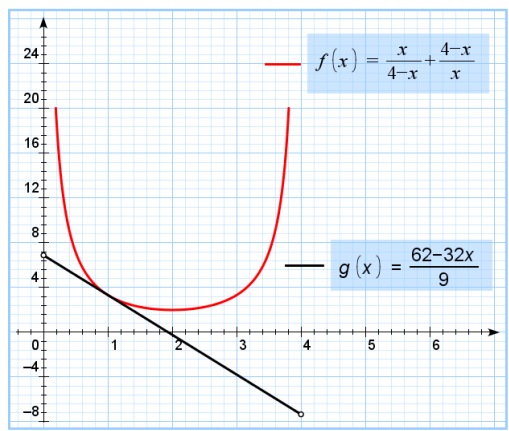

Fig. 8

Indeed, multiply the left and right-hand sides of (8) by a positive 9x(4-x), then after collecting similar terms, we obtain

$$2x^3 - 13x^2 + 20x - 9 \leq 0,$$

which factorizes as

$$(x-1)^2(2x-9) \leq 0. \tag{9}$$

As *x*<4, then 2*x*<8<9, and (9) is true, then inequality (8) holds. Applying (8), we obtain

$$f(a) + f(b) + f(c) + f(d) \geq \frac{62 \cdot 4 - 32(a+b+c+d)}{9} = \frac{40}{3},$$

and equality is achieved when *a*=*b*=*c*=*d*=1. Consequently, the lowest value of the expression

$$\frac{a}{b+c+d} + \frac{b}{c+d+a} + \frac{c}{d+a+b} + \frac{d}{a+b+c} + \frac{b+c+d}{a} + \frac{c+d+a}{b} + \frac{d+a+b}{c} + \frac{a+b+c}{d},$$

where *a*, *b*, *c*, *d* are positive numbers summing to $\frac{40}{3}$.

As previously mentioned, the choice of the tangency type is defined by the tangent inequality to be proved. One of the clearest examples demonstrating this instance is the following problem.

*Example* 9 (Chetkovski, 2012: **202**, problem 185). Let *a*, *b*, *c* >0 be real numbers such that $a^{2/3} + b^{2/3} + c^{2/3} = 3$. Prove the inequality

$$a^2 + b^2 + c^2 \geq a^{4/3} + b^{4/3} + c^{4/3}.$$

*Solution.* Let $f(x) = x^2 - x^{4/3}$, $g(x) = kx^{2/3} + m$, where $0 < x < 3\sqrt{3}$. Note that at *a*=*b*=*c*=1 our inequality becomes an equality. Numbers *k* and *m* are chosen such that $f(1) = g(1)$, $f'(1) = g'(1)$. That is, $0 = k + m$, $\frac{2}{3} = \frac{2}{3}k$. Consequently, $g(x) = x^{2/3} - 1$. As the inequality $x^2 - x^{4/3} \geq x^{2/3} - 1$ (see Figure 9) is equivalent to $(x^{2/3} - 1)^2(x^{2/3} + 1) \geq 0$, then

$$a^2 + b^2 + c^2 - a^{4/3} - b^{4/3} - c^{4/3} \geq a^{2/3} + b^{2/3} + c^{2/3} - 3 = 0.$$

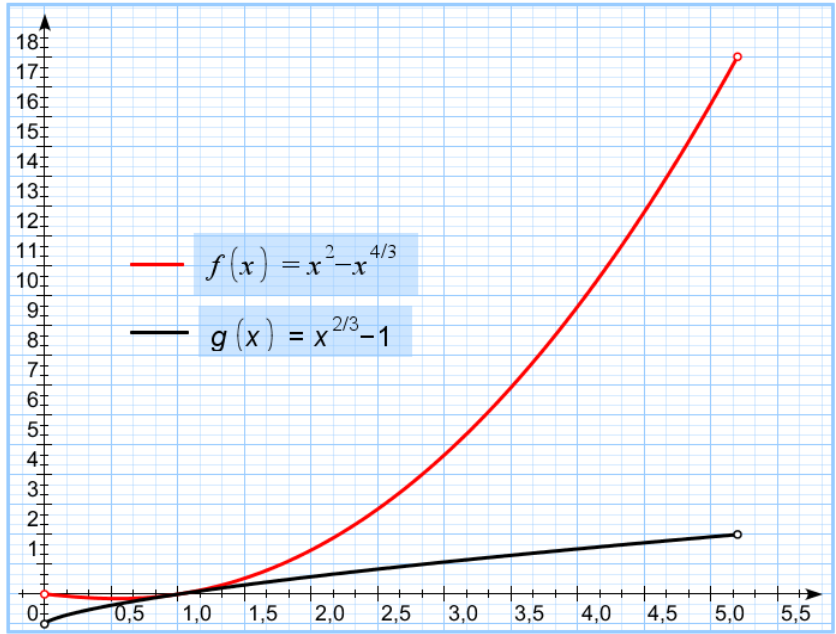

Fig. 9

*Example* 10 (Pham Kim Hung, 2007: **37**, ex. 2.1.5). Let $a$, $b$, $c$, $d$ be non-negative real numbers. Prove that

$$\frac{a}{b^2+c^2+d^2} + \frac{b}{a^2+c^2+d^2} + \frac{c}{a^2+b^2+d^2} + \frac{d}{a^2+b^2+c^2} \geq \frac{4}{a+b+c+d}.$$

*Solution.* According to the Cauchy-Schwarz inequality, we have

$$\left(\frac{a}{b^2+c^2+d^2} + \frac{b}{a^2+c^2+d^2} + \frac{c}{a^2+b^2+d^2} + \frac{d}{a^2+b^2+c^2}\right)(a+b+c+d) \geq$$

$$\geq \left(\sqrt{\frac{a^2}{b^2+c^2+d^2}} + \sqrt{\frac{b^2}{a^2+c^2+d^2}} + \sqrt{\frac{c^2}{a^2+b^2+d^2}} + \sqrt{\frac{d^2}{a^2+b^2+c^2}}\right)^2.$$

Therefore it is sufficient to prove the inequality

$$\sqrt{\frac{a^2}{b^2+c^2+d^2}} + \sqrt{\frac{b^2}{a^2+c^2+d^2}} + \sqrt{\frac{c^2}{a^2+b^2+d^2}} + \sqrt{\frac{d^2}{a^2+b^2+c^2}} \geq 2. \tag{10}$$

Assume that $a^2+b^2+c^2+d^2=4$. In this case, it can be written as follows

$$\frac{a}{\sqrt{4-a^2}} + \frac{b}{\sqrt{4-b^2}} + \frac{c}{\sqrt{4-c^2}} + \frac{d}{\sqrt{4-d^2}} \geq 2. \tag{11}$$

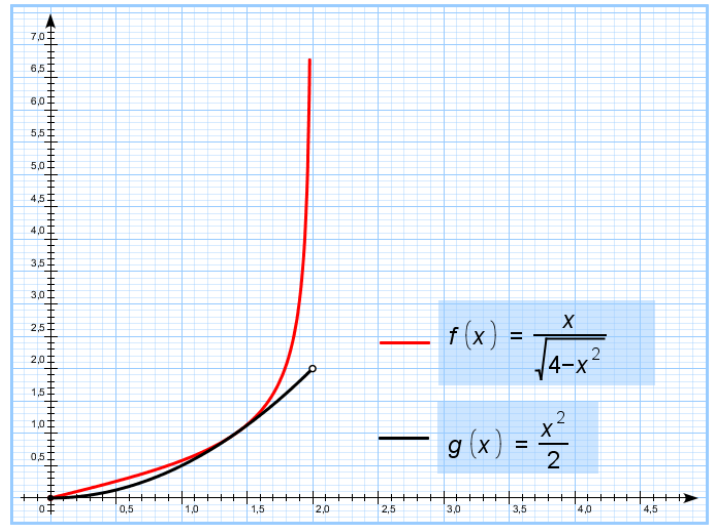

Fig. 10

Let $f(x) = \frac{x}{\sqrt{4-x^2}}$, $g(x)=kx^2+m$, where $0 \leq x < 2$. Note that the inequality (11) becomes an equality for $a = b = \sqrt{2}$, $c=d=0$. Therefore, the numbers $k$ and $m$ are chosen so that $f(0) = g(0)$, $f(\sqrt{2}) = g(\sqrt{2})$, $f'(\sqrt{2}) = g'(\sqrt{2})$. It is not difficult to make sure that the specified conditions are satisfied only by function $g(x) = \frac{x^2}{2}$, and for every $x \in [0, 2)$ we have $\frac{x}{\sqrt{4-x^2}} \geq \frac{x^2}{2}$ (see Figure 10), as it is equivalent to $x^2(x^2 - 2)^2 \geq 0$. Consequently,

$$\frac{a}{\sqrt{4-a^2}} + \frac{b}{\sqrt{4-b^2}} + \frac{c}{\sqrt{4-c^2}} + \frac{d}{\sqrt{4-d^2}} \geq \frac{a^2+b^2+c^2+d^2}{2} = 2.$$

QED■

**References**

[1] Chetkovski, Z. (2012) *Inequalities. Theorems, Techniques and Problems,* Berlin Heidelberg: Springer Verlag.
[2]. Li K. *Using tangent lines to prove inequalities*. Mathematical Excalibur. 2005-2006. V.10. No. 5. P.1.
[3] Pham Kim Hung *Secrets in Inequalities* (*Volume* 1). – Zalău : Gill. 2007.
[4] Xiong Bin and Lee Peng Yee (2007) *Mathematical Olympiad in China Problems and Solutions*, Shanghai: Word Scientific.